\documentclass[10pt]{article}

\usepackage{amsmath,amscd,amssymb,amsfonts,graphicx,color}
\usepackage{verbatim}
\usepackage[affil-it]{authblk}

\def\R     {\ensuremath{\mathbb R}}
\def\N     {\ensuremath{\mathbb N}}

\begin{document}

\date{\today}
\title{Matrix methods for radial Schr\"{o}dinger eigenproblems defined on a
semi-infinite domain}

\author[a]{Lidia Aceto}   
\affil[a]{Dipartimento di Matematica, Universit\`a di Pisa, Italy} 
\author[a]{Cecilia Magherini}
 \author[b]{Ewa\,B.\ Weinm\"{u}ller}    
\affil[b]{Department for Analysis and Scientific Computing,  Vienna University
of Technology, Austria}

\maketitle

\begin{abstract}
In this paper, we discuss numerical approximation of the eigenvalues of the
one-dimensional radial Schr\"{o}dinger equation posed on a semi-infinite interval.
The original problem is first transformed to one defined on a finite domain by
applying suitable change of the independent variable. The eigenvalue problem for the
resulting differential operator is then approximated by a generalized algebraic eigenvalue problem
arising after discretization of the analytical problem by the matrix method
based on high order finite difference schemes. Numerical experiments illustrate
the performance of the approach.

\bigskip


\smallskip

\end{abstract}

\section{Introduction}
The aim of this paper is to investigate certain aspects arising in the numerical treatment
of the following eigenvalue problem (EVP):
\begin{eqnarray}
&& - u''(r)+\left(  \frac{\ell (\ell+1)}{r^2} + V(r) \right) u(r) = \lambda
u(r), \quad r \in (0,\infty), \label{prob}
\end{eqnarray}
subject to boundary conditions
\begin{equation}
u(0) = u(\infty)=0, \label{BC}
\end{equation}
where $\ell \in \N$, the function $V(r)$ satisfies $\displaystyle{\lim_{r \to
\infty} V(r)=0},$ $\lambda$ is an eigenvalue, and $u(r)$ is the associated
eigenfunction.

Equation (\ref{prob}) is known in the literature as {\em radial} Schr\"{o}dinger equation
with {\em underlying potential} $V(r).$ An important example of this type of problems is
the hydrogen atom equation corresponding to $V(r) = -Z/r$ with $Z >0.$\\

Many currently available numerical techniques to handle this problem are based on the so-called
{\em regularization} which, in our context, may mean replacing (\ref{BC}) by
\begin{equation}\label{BCt}
u(\varepsilon) = u(R) = 0,
\end{equation}
where $\varepsilon$ is strictly positive and small and $R$ is large.
Clearly, equation (\ref{prob}) subject to (\ref{BCt}) is
a regular Sturm-Liouville problem on $[\varepsilon,R]$ and classical methods can be
used to approximate its eigenvalues. The accuracy of these approximations,
however, strongly depends on the choice of the
cutoff points $\varepsilon$ and $R$ as discussed, for example, in \cite{pryce}.
Concerning the choice of $R,$ a generalization of the so-called WKB-approximation
introduced in \cite{ixa} for nonharmonic oscillators, was proposed in \cite{lirvv}.
In case of a problem whose potential has a Coulomb-like tail, the authors proposed
to impose suitably adapted boundary conditions at the right endpoint $R$ which allowed a
noticeable reduction of the size of $R.$ On the other hand, an ad-hoc procedure,
treating $V(r)$ as a perturbation of the reference potential $\ell (\ell+1) /r^2$ in a
neighborhood of the origin, was investigated in \cite{imv,lirvv}.
The aim of this procedure was to find an approximation for $u(\varepsilon)$ and $u'(\varepsilon).$
All these estimates were then used in a two-sided shooting procedure.\\

The approach proposed here is different. We first apply suitable change of variable, $t(r)$,
in order to transform the problem posed on $r \in (0, \infty)$ to a problem
posed on a finite domain $t \in (0,1).$ After the transformation the EVP assumes the following general form:
\begin{eqnarray} \label{sisewa0}
- A_2(t) \, {\mathbf{v}}'' (t) + A_1(t) \,  {\mathbf{v}}' (t)+ A_0(t) \,  \mathbf{v}(t) = \lambda \, \mathbf{v}(t), \quad t \in (0,1),
\end{eqnarray}
where  $A_i(t), \, i=0,1,2,$  are singular at $t=0$ and/or $t=1.$ The EVP (\ref{sisewa0}) is then
augmented by suitable boundary conditions. We use results provided in
\cite{deHoogWeiss79,deHoogWeiss80b,weinm84} to show that the above singular EVP is
well-posed and to describe the smoothness of its solution.

Numerical approximations of the eigenvalues are obtained by applying the so-called matrix
methods which transform the EVP for the differential operator into a generalized algebraic EVP.
More precisely, the equation (\ref{sisewa0}) is discretized in its original second order formulation
by using the finite difference schemes introduced in \cite{as}. We stress that
the application of these methods is possible in spite of the singularities at
$t=0,1$, since the corresponding discrete problem does
not involve the values of
the coefficients functions in (\ref{sisewa0}) at the interval endpoints.\\

It is worth mentioning that the idea of reducing the continuous problem to a
finite domain was already used in the development of the codes SLEIGN2 \cite{bez} and SLF02F \cite{mp}.
Nevertheless, the numerical schemes used in the implementation of the shooting procedure make
use of the coefficient functions at the endpoints. This means that, in our case, where such
functions become unbounded, cutting off the interval ends becomes inevitable. \\

We have organized the paper as follows. In Section \ref{sec2}, we propose two ways of changing
the independent variable for the transformation of the original problem to a finite domain and
discuss the properties of the resulting singular BVPs and EVPs. In Section \ref{sec3}, we describe
in some detail the numerical procedure based on the matrix method. Finally,
Section \ref{sec4} contains 
the results of the numerical simulation for the hydrogen atom equation and
models studied in \cite{Roy}. Here, we also show numerical results related to a
third change of independent variable, which is analyzed in detail in
\cite{acetoetal2013}.

\section{Reformulation of the problem on a finite domain}\label{sec2}
The first question we would like to address is how to transform problem (\ref{prob})-(\ref{BC})
posed on the semi--infinite interval to a finite domain.
In general, if $t(r)$ is given and  $u(r)=:z(t(r)),$ then we can rewrite
(\ref{prob}) as
\begin{eqnarray} \label{probmod}
- \frac{d^2}{dt^2}  z(t) \left( \frac{d}{dr} t(r)  \right)^2 \!\!-\! \frac{d}{dt}  z(t) \frac{d^2}{dr^2} t(r) \!+\!
\left(  \frac{\ell (\ell+1)}{r^2} \!+\! V(r) \right)  z(t) \!=\! \lambda z(t).
\end{eqnarray}

\subsection{Transformation doubling the size of the ODE system: TDS}
The transformation TDS is based on the following change of the independent variable:
\begin{eqnarray}
\label{cvewa}
&&t(r)=\frac{1}{r}, \quad r \in[1,\infty).
\end{eqnarray}
We use (\ref{cvewa}) to reformulate (\ref{probmod}) as follows:
$$
 - z''(t) -  \frac{2}{t} z'(t)  +
\left( \frac{ \ell (\ell+1)}{ t^2} + \frac{1}{t^4}V \left(\frac{1}{t}\right) \right)  z(t) = \lambda  \frac{ z(t)}{t^4}, \quad t \in (0,1],
$$
and therefore (\ref{prob}) posed on the interval $(0,\infty)$ can be
transformed to the finite interval,
\begin{eqnarray} \nonumber
\begin{array}{l}
\displaystyle{
-  u''(t)+\left(  \frac{\ell (\ell+1)}{t^2} + V(t) \right) u(t) = \lambda u(t),
\quad t \in (0,1],} \\
\\
\displaystyle{
 - z''(t) -  \frac{2}{t} z'(t)  +
\left( \frac{ \ell (\ell+1)}{ t^2} + \frac{1}{t^4}V \left(\frac{1}{t}\right)
\right)  z(t) = \lambda  \frac{ z(t)}{t^4}, \quad t \in (0,1].}
\end{array}
\end{eqnarray}
In matrix notation, this system of equations can be written as
\begin{eqnarray} \label{sisewa}
- {\mathbf{v}}'' (t) + {\tilde A}_1(t) \,  {\mathbf{v}}' (t)+ {\tilde A}_0(t)
\, \mathbf{v}(t) = \lambda \, B (t) \, \mathbf{v}(t), \quad t \in (0,1],
\end{eqnarray}
with $\mathbf{v}(t) = (u(t), z(t))^T$ and
\begin{eqnarray*}
{\tilde A}_1(t) &=&   \left(\begin{array}{cc}
 0 & 0\\
 0 & - 2 t^{-1}
\end{array}\right),  \\
{\tilde A}_0(t) &=&   \left(\begin{array}{cc}
\ell (\ell+1) t^{-2}+V(t)& 0\\
 0 & \ell (\ell+1) t^{-2} +  t^{-4}V \left(t^{-1}\right)
\end{array}\right), \\
B(t) &=&  \left(\begin{array}{cc}
 1 & 0\\
 0 &  t^{-4}
\end{array}\right).
\end{eqnarray*}

Note that $B(t)$ is nonsingular for $t \in (0,1),$ and hence, (\ref{sisewa}) can be written
in the general form (\ref{sisewa0}).\\

In the sequel, we investigate if the above singular EVP is well-posed.
This is done by first examining the boundary conditions. To this aim, we follow the arguments
from \cite{deHoogWeiss76,deHoogWeiss80b}. Although, we will numerically simulate
the EVPs in form (\ref{sisewa0}), for the analysis, we have to rewrite the problem into its first order form.
It turns out that here $t=0$ is a singular point and therefore, we have to
investigate the local behavior of the ODE in the vicinity of this point.

If we transform (\ref{sisewa}) to a first order system
of ODEs for the new vector
\begin{equation}\label{vecyTDS}
\mathbf{y}(t)=(y_1(t),y_2(t),y_3(t),y_4(t))^T:=(\mathbf{v}(t), t
\mathbf{v'}(t))^T \in \R^4,
\end{equation}
then we obtain
\begin{eqnarray} \label{sisewa4}
t^4 {\mathbf{y'}}(t)  - M(t){\mathbf{y}}(t) = \lambda \, G(t) \, \mathbf{y}(t), \quad t \in (0,1],
\end{eqnarray}
where the matrices $M(t)$ and $G(t)$ include the data from (\ref{sisewa}), namely
\begin{eqnarray*}
&&M(t) =   \left(\begin{array}{cccc}
0 & 0 & t^3 & 0 \\
0 & 0 & 0 & t^3 \\
t^3\ell (\ell+1) + t^5 V(t) & 0 & t^3 & 0 \\
 0 & t^3 \ell (\ell+1) + t V \left( t^{-1} \right) & 0 & -t^3
\end{array}\right), \\
&&G(t) =  \left(\begin{array}{cccc}
0 & 0 & 0 & 0 \\
0 & 0 & 0 & 0 \\
-t^5 & 0 & 0 & 0 \\
 0 & -t  & 0 & 0
\end{array}\right).
\end{eqnarray*}
For the investigation of the local behavior of (\ref{sisewa4}) around $t=0$, note that
${\displaystyle \lim_{t \to 0^+} V \left(t^{-1} \right)= \lim_{r \to
\infty} V(r) = 0}$.
Moreover, we assume that $\displaystyle{\lim_{t \to 0^+} t^5 V(t)=0,}$ and the
higher derivatives of $t^5 V(t)$ exist
and are continuous on $[0,1]$. This means that in (\ref{sisewa4}), $M(t) = M + A(t)$ and $G(t) = N + C(t)$,
where $M=M(0)$ and $N=G(0)$
are zero matrices. Consequently, (\ref{sisewa4}) has the form
\begin{eqnarray}
t^4 {\mathbf{y'}}(t)  - A(t){\mathbf{y}}(t) = \lambda \, C(t) \, \mathbf{y}(t), \quad t \in (0,1]. \label{sisewa5}
\end{eqnarray}
The associated boundary conditions read:
\begin{equation}
y_1(0)=y_2(0)=0, \quad y_1(1)=y_2(1), \quad y_3(1)=-y_4(1), \label{bcsisewa6}
\end{equation}
which is equivalent to
\begin{equation}
\!\!\!B_0\mathbf{y}(0)+B_1\mathbf{y}(1)=\mathbf{0}, 
B_0=\left(\begin{array}{cccc}
1 & 0 & 0 & 0 \\
0 & 1 & 0 & 0 \\
0 & 0 & 0 & 0 \\
0 & 0 & 0 & 0
\end{array}\right),  B_1=\left(\begin{array}{rrrr}
0 & 0 & 0 & 0 \\
0 & 0 & 0 & 0 \\
1 & -1 & 0 & 0 \\
0 & 0  & 1 & 1
\end{array}\right). \label{bcsisewa7}
\end{equation}
First note that the form of the EVP (\ref{sisewa5})-(\ref{bcsisewa6}) corresponds exactly to the one
of the EVP (1.2) studied in \cite{deHoogWeiss80b}. To discuss the boundary conditions, we have to
look at the associated BVP
\begin{eqnarray}
t^4 {\mathbf{y'}}(t)  - M(t){\mathbf{y}}(t)=t^4 {\mathbf{y'}}(t)  -
    (M + A(t)){\mathbf{y}}(t) = {\mathbf{g}}(t), \quad t \in (0,1], \label{sisewa6}
\end{eqnarray}
subject to (\ref{bcsisewa7}), cf.  \cite[problem (3.9)]{deHoogWeiss80b}. Since
the matrix $M=M(0)$ is a zero matrix,
its eigenvalues are $\mu_1=\mu_2=\mu_3=\mu_4:=\mu=0$ and the corresponding
eigenspace is $\R^4$.
Consequently, the orthogonal projection $R$ onto the eigenspace of $M$ associated with $\mu=0$ is $R=I$ and
the condition $(I-R)A(0)=O$ is satisfied, see \cite[requirement
(3.10)]{deHoogWeiss80b}. Moreover, due to \cite[Theorem 3.2]{deHoogWeiss80b},
for any ${\mathbf{g}}, A \in C[0,1]$  there exists a unique solution
${\mathbf{y}} \in C^1[0,1]$ 
of the BVP (\ref{bcsisewa7})-(\ref{sisewa6}) since
$\mbox{rank}[B_0, B_1] := k = 4$, and the linear differential operator 
$t^4 {\mathbf{y'}}(t)  - A(t){\mathbf{y}}(t)$ is Fredholm with index equal to
$\mbox{rank}[R] - k =0$. This result immediately carries over to the
EVP problem (\ref{sisewa5})-(\ref{bcsisewa6}). Here, $(I-R)C(0) = O $ holds,
cf.  \cite[(7.1)]{deHoogWeiss80b} and $A, C$ are smooth functions. Therefore,
according to \cite[Theorem 7.2]{deHoogWeiss80b} the EVP is well-posed and has a
solution in $C^\infty[0,1]$. \\

For the numerical treatment, we use system (\ref{sisewa}), where the second
equation is premultiplied by $t^4,$ together with boundary conditions, see
(\ref{vecyTDS}) and  (\ref{bcsisewa6}),
\begin{equation}
\mathbf{v}(0)=\mathbf{0}, \qquad (1,-1)\,\mathbf{v}(1)=0, \quad (1,1)\,\mathbf{v}'(1)=0. \label{bcsisewa5}
\end{equation}

\subsection{Transformation compressing the infinite interval: TCII}
We now consider an alternative change of independent variable described by
\begin{eqnarray} \label{cvnoi}
&&t(r)=\frac{r}{r+\xi},  \quad \xi>0, \quad r \in [0,\infty).
\end{eqnarray}
Using (\ref{cvnoi}) in (\ref{probmod}) yields the following new form of (\ref{prob}):
\begin{equation} \label{probmodnoi}
\!\!- z''(t) + \frac{2}{1-t}   z'(t)+
\left( \frac{ \ell (\ell+1)}{ t^2 (1-t)^2} + \frac{\xi^2}{(1-t)^4}V \left( \frac{\xi t}{1-t}\right) \right)  z(t)
= \frac{ \lambda  \xi^2}{(1-t)^4} z(t),
\end{equation}
subject to boundary conditions
\begin{equation}  \label{bcnoi} z(0)=z(1)=0.
\end{equation}
First of all we note that there are two critical points $t=0$ and $t=1$
in the differential operator in (\ref{probmodnoi}). Our aim is to show that boundary
conditions (\ref{bcnoi}) are posed in such a way that the associated BVP is well-posed. To this aim,
we have to investigate the ODE in the vicinity of
$t=0$ and $t=1$. Let us first consider $t=0$.
Setting
$$
a_1(t) = \frac{2t}{1-t}, \qquad a_0(t) =  \frac{ \ell (\ell+1)}{(1-t)^2} \!+\! \frac{t^2\xi^2}{(1-t)^4}V \left( \frac{\xi t}{1-t} \right),
\qquad b(t)= - \frac{\xi^2}{(1-t)^4},
$$
we rewrite (\ref{probmodnoi}) to obtain
the form
$$
 z''(t) - \frac{a_1(t)}{t} z'(t) - \frac{a_0(t)}{t^2} z(t) = \lambda b(t) z(t),
$$
and transform it to the following first order system for the vector
$\mathbf{y}(t)=(y_1(t),y_2(t))^T:=(z(t),tz'(t))^T$
\begin{equation} \label{sysevp2}
t \mathbf{y}'(t) - M(t)\mathbf{y}(t) = \lambda G(t) \mathbf{y}(t),
\end{equation}
subject to
\begin{equation}
B_0{\mathbf{y}}(0)+B_1{\mathbf{y}}(1)={\bf 0}, \qquad
B_0=\left(\begin{array}{cc}
1 & 0  \\
0 & 0
\end{array}\right), \quad B_1=\left(\begin{array}{cc}
0 & 0 \\
1 & 0
\end{array}\right), \label{bcsisewa8}
\end{equation}
where
\begin{equation} \label{data}
M(t) = \left( \begin{array}{cc}
                 0 & 1 \\
         a_0(t) & 1+a_1(t)
               \end{array} \right), \qquad
G(t) = \left( \begin{array}{cc}
                 0 & 0 \\
                t^2 b(t) & 0
               \end{array} \right).
\end{equation}
Also here, if we assume that $\displaystyle{\lim_{r\to 0^+} rV(r)}$ is finite
we have $M(t) = M + A(t)$ and $G(t) = N + C(t)$, where $N$ is a zero
matrix
and
$$
M=M(0) = \left( \begin{array}{cc}
                 0 & 1 \\
         a_0(0) & 1+a_1(0)
               \end{array} \right) =
        \left( \begin{array}{cc}
                 0 & 1 \\
         \ell(\ell+1) & 1
               \end{array} \right).
$$

In contrast to (\ref{sisewa4}), where due to $t^4$ in the leading term the ODE
admits a singularity of the second kind, in (\ref{sysevp2}) a singularity of the
first kind arises. Therefore, we can apply results from \cite{deHoogWeiss76} to
analyze the boundary conditions of the
problem. We first calculate the eigenvalues of the matrix $M(0)$ and obtain 
$\mu_1 = -\ell \leq 0$ and  $\mu_2=1+\ell > 0.$ Let us focus on two cases used
in the numerical simulations.
\begin{enumerate}
\item[Case 1: $\ell=0$] For $\ell=0$, the eigenvalues of the matrix $M$ are
$\mu_1=0$ and $\mu_2=1.$
First, we have to calculate the related eigenvectors ${\mathbf{w_1}}$ and
${\mathbf{w_2}}$ and construct
two projection matrices  $R$ and $S,$ $R+S=I_2$, projecting onto eigenspaces
of $M$ associated with
$\mu_1$ and $\mu_2$, respectively. This yields
$$
{\mathbf{w_1}} = \left( \begin{array}{r}
                 1 \\
                 0   \end{array} \right), 
{\mathbf{w_2}} = \left( \begin{array}{r}
                 1 \\
                 1
               \end{array} \right), 
R = \left( \begin{array}{rr}
                 1 & -1 \\
                 0 & 0
               \end{array} \right),
S = \left( \begin{array}{rr}
                 0 & 1 \\
                 0 & 1
               \end{array} \right).
$$
According to  \cite[Theorem 3.2]{deHoogWeiss76}, the linear operator $t \mathbf{y}'(t) - M(t)\mathbf{y}(t)$
is Fredholm with index equal to $\mbox{rank}(R+S)-\mbox{rank}[B_0R,B_1]=2-2=0$
since
\[
{\mbox{rank}}[B_0R,B_1]={\mbox{rank}}  \left[  \left(\begin{array}{rrrr}
1 & -1  & 0 & 0 \\
0 & 0   & 1 & 0
\end{array}    \right) \right] =2.
\]
Again, this means that for any $\mathbf{g}, A \in C[0,1]$ the BVP,
\begin{equation} \label{sysevp41}
t \mathbf{y}'(t) - M(t)\mathbf{y}(t) = \mathbf{g}(t), \quad t \in (0,1), \quad
B_0{\mathbf{y}}(0)+B_1{\mathbf{y}}(1)={\bf 0},
\end{equation}
where the problem data has been specified in (\ref{bcsisewa8})-(\ref{data}),
is well-posed and has as solution ${\mathbf{y}} \in C[0,1] \cap C^1(0,1]$. We have an analogous
result for the EVP (\ref{sysevp2})-(\ref{bcsisewa8}) with $A, G \in C[0,1]$.
Since the positive eigenvalue of $M$ is relatively small, we would need further investigations to
show that also higher derivatives of $\mathbf{y}$ are smooth, cf. \cite[Theorem
10.2]{deHoogWeiss76}.
\item[Case 2: $\ell=3$] Here, the eigenvalues of the matrix $M$ are $\mu_1=-3$ and $\mu_2=4$.
Again, we first calculate the related eigenvectors ${\mathbf{w_1}}$ and
${\mathbf{w_2}}$ and construct
two projection matrices $Q$ and $S$, $Q+S=I_2$, projecting onto eigenspaces of $M$ associated with
$\mu_1$ and $\mu_2$, respectively. This yields
\[
{\mathbf{w_1}} = \left( \begin{array}{r}
                 1 \\
                -3   \end{array} \right), \quad
{\mathbf{w_2}} = \left( \begin{array}{r}
                 1 \\
                 4
               \end{array} \right), \quad
Q = \frac{1}{7}\left( \begin{array}{rr}
                 4 & -1 \\
                 -12 & 3
               \end{array} \right), \quad
S = \frac{1}{7}\left( \begin{array}{rr}
                 3 & 1 \\
                 12 & 4
               \end{array} \right).
\]
First of all, $S{\mathbf{y}} \in C^1[0,1]$ and $S{\mathbf{y}}(0)= S{\mathbf{y}}'(0)=\mathbf{0}$,
see \cite[Lemma 3.5]{deHoogWeiss76}. Moreover, condition $Q\mathbf{y}(0)=\mathbf{0}$ is necessary and
sufficient for $\mathbf{y}$ to be in $C[0,1]$. To see that this condition is satisfied,
we have to take into account that the ODE in (\ref{sysevp41}) arises from
$$
z''(t) - \frac{a_1(t)}{t} z'(t) - \frac{a_0(t)}{t^2} z(t) = g(t), \quad t \in (0,1],
$$
and thus $\mathbf{g}(t) = t^2(0,g(t))^T$. Using the special structure of
$\mathbf{g}$ and \cite[Lemma 3.1]{weinm84}, we see that from $z(0)=y_1(0)=0$,
$y_1 \in C^1[0,1]$ follows and therefore
$$
(1\,,\, 0)\,Q\,\mathbf{y}(0) = \frac{4}{7}y_1(0) - \frac{1}{7}y_2(0) =
-\frac{1}{7} \lim_{t \to 0^+} ty_1'(t) = 0
$$
holds.  \\

Now, according to \cite[Theorem 3.2]{deHoogWeiss76}, the linear operator $t \mathbf{y}'(t) - M(t)\mathbf{y}(t)$
is Fredholm with index equal to $\mbox{rank}(S)-\mbox{rank}[B_0R,B_1]=1-1=0$
since the orthogonal projection $R$ onto the eigenspace of $M$ associated with $\mu =0 $ is zero  and
\[
{\mbox{rank}}[B_0R,B_1]= {\mbox{rank}}[B_1] = 
{\mbox{rank}}\left[\left(\begin{array}{rr}
0 & 0 \\
1 & 0
\end{array}\right)\right] =1.
\]
Thus, for any $\mathbf{g}, A \in C^3[0,1]$ the BVP (\ref{sysevp41})
with the problem data given in (\ref{bcsisewa8})-(\ref{data}) is well-posed
and has as solution ${\mathbf{y}} \in C^3[0,1] \cap C^4(0,1]$.
We have an analogous result for the EVP (\ref{sysevp2})-(\ref{bcsisewa8})
for $A, G \in C^3[0,1]$. Since the positive eigenvalue of $M$ is slightly larger than in Case 1,
we can show more smoothness in $\mathbf{y}$, cf.  \cite[Theorem 10.2]{deHoogWeiss76}.
\end{enumerate}
Similar investigations for $t=1$ show that this point is not a critical point and the solution is analytic at $t=1$,
see \cite{acetoetal2013,amodioetal2013}. \\

For the numerical experiments, we use (\ref{probmodnoi})
premultiplied by $(1-t)^4/\xi^2$ together with boundary conditions (\ref{bcnoi}).

\section{Finite difference schemes} \label{sec3}
The numerical methods that we have used discretize equation (\ref{sisewa0})
in its original second order formulation. In particular, given the uniform mesh
$$
t_i=i \,h, \quad  i=0,1,\dots,N+1,\quad h=1/(N+1),
$$
for the interval $[0,1],$ the first and second order derivatives of the solution
at the inner grid points are approximated by applying suitable
$(2k)$-step finite difference schemes introduced in \cite{as}. More precisely, for each
$i=k,k+1,\dots,N+1-k,$
\begin{equation} \label{v1-v2}
 {\mathbf{v}}' (t_i) \approx    \frac{1}{h} \sum_{j=-k}^k
\beta_{j+k}{\mathbf{v}}_{i+j},   \quad
{\mathbf{v}}'' (t_i) \approx \frac{1}{h^2}
\sum_{j=-k}^k \gamma_{j+k} {\mathbf{v}}_{i+j},
\end{equation}
where ${\mathbf{v}}_i \approx {\mathbf{v}}(t_i),$ for each $i.$ The
coefficients $\{\beta_j\}_{j=0}^{2k}$ and $\{\gamma_j\}_{j=0}^{2k}$ are uniquely
determined by imposing the formulas to be of consistency order $2k.$  The resulting methods
turn out to be symmetric, i.e.,  $\beta_j = - \beta_{2k-j}$ and $\gamma_j =  \gamma_{2k-j},$
for each $j=0,1,\dots,k.$ In particular, the $2$-step
schemes coincide with the ones used in \cite[Section~5.3]{rep}.
Using the terminology of Boundary Value Methods the formulas in (\ref{v1-v2})
are called {\em main methods} \cite{BT}.
When $k>1,$ these formulas are augmented by suitable {\em initial} and {\em final additional methods}
which provide approximations of the first and second order derivatives at the meshpoints close to the interval ends.
In particular, for each  $i=1,2,\dots,k-1,$
$$
 {\mathbf{v}}' (t_i) \approx  \frac{1}{h}
\sum_{j=0}^{2k}
\beta_j^{(i)}{\mathbf{v}}_j, \quad
{\mathbf{v}}'' (t_i) \approx \frac{1}{h^2}
\sum_{j=0}^{2k+1}
\gamma_j^{(i)}{\mathbf{v}}_j,
$$
while, for each  $i=N+2-k,\dots,N$ and $r=N+1-2k,$
$$
{\mathbf{v}}' (t_i) \approx \frac{1}{h}
\sum_{j=0}^{2k}
\beta_j^{(i-r)}{\mathbf{v}}_{r+j}, \quad
{\mathbf{v}}'' (t_i) \approx   \frac{1}{h^2}
\sum_{j=0}^{2k+1}
\gamma_j^{(i-r)}{\mathbf{v}}_{r+j-1}.
$$
The involved coefficients are determined by requiring that the additional schemes are of
the same order as the main formulas, i.e. $2k,$ \cite{as}. It is worth
mentioning that such schemes have been already used for solving singular
Sturm-Liouville
problems in \cite{agm,as0}.\\
Since for both transformations, TDS and TCII, $\mathbf{v}(0) = \mathbf{0}$
holds, the following
system of equations arises
after the discretization of (\ref{sisewa0}):
\begin{eqnarray} \label{eqdis}
\mathcal{R}  \hat{\mathbf{v}} := \left( -D_2 \left(  \hat{\Gamma} \otimes  I_m \right)  +D_1 \left( \hat{B}
\otimes I_m \right)  + \hat{D}_0 \right) \hat{\mathbf{v}} = \lambda \mathbf{v}.
\end{eqnarray}
Here $I_m$ is the identity matrix of dimension $m,$ with $m=2,1$ for TDS and TCII, respectively,
\begin{eqnarray}
D_i &=& \mbox{blockdiag}\left( A_i(t_1),A_i(t_2),\dots, A_i(t_N) \right), \quad i=0,1,2, \nonumber \\
\hat{D}_0 &=& \left(D_0 ~\,|~~ {\bf 0}_N \otimes I_m \right), \quad {\bf 0}_N =
\left(0,\ldots,0\right)^T \in \R^N,\nonumber \\
   \hat{\mathbf{v}}^T &=& \left(\mathbf{v}^T \,\,  \mathbf{v}_{N+1}^T \right)
  = \left(\mathbf{v}_1^T \, \dots \,\mathbf{v}_N^T \,  \mathbf{v}_{N+1}^T \right) \label{vetv}
  \approx   \left( \mathbf{v}^T(t_1) \, \dots \, \mathbf{v}^T(t_N)\, \mathbf{v}^T (t_{N+1}) \right).
\end{eqnarray}

\noindent Finally, $\hat{\Gamma}, \, \hat{B} \in \R^{N \times (N+1)}$ contain the coefficients of the
difference schemes.
For example, for the method of order $4,$
\begin{eqnarray*}
 \hat{\Gamma} &=& {
 \frac{1}{12 \, h^2} \left(\begin{array}{rrrrrrr}
  - 15&  - 4 &  14 & -6 & 1 \\
  16&  - 30 &  16 & - 1 \\
 -1& 16&  - 30 &  16 & - 1 \\
 &\ddots& \ddots& \ddots& \ddots& \ddots \\
&&- 1&  16&  - 30 &  16 & - 1 \\
& 1 &  -6 &  14 &  -4 & -15 & 10
               \end{array} \right),} \\
\hat{B} &=&  { \frac{1}{12 \, h} \left(\begin{array}{rrrrrrr}
  -10 & 18 & -6 & 1  \\
  - 8 & 0 &  8 & -1  \\
 1& - 8 & 0 &  8 & -1  \\
&\ddots& \ddots& \ddots& \ddots &\ddots  \\
&&  1&  - 8 & 0 &  8 & -1  \\
 && -1&  6& -18 & 10 & 3
             \end{array} \right).}
\end{eqnarray*}
Let us now describe the discretization of the boundary conditions at $t=1.$
We have to distinguish between TDS and TCII. For TDS, the last two conditions in (\ref{bcsisewa5})
are approximated as follows:
\begin{equation} \label{bcewadis}
 \left(1, \, -1 \right) \mathbf{v}_{N+1} =0, \qquad
 \left(1, \, 1 \right) \mathbf{v}'(1) \approx
\left(1, \, 1 \right)\sum_{j=0}^{2k} \beta_j^{(2k)}    \mathbf{v}_{N-2k+j}  = 0,
\end{equation}
where $\beta_j^{(2k)}$ are the coefficients of the classical $(2k)$-step BDF method.
Equation (\ref{eqdis}) augmented with (\ref{bcewadis}) form the following
generalized algebraic EVP,
$$
  \mathcal{R}_1 \,  \hat{\mathbf{v}} = \lambda \mathcal{S} \,   \hat{\mathbf{v}}, \qquad \mathcal{S} =   \left( \begin{array}{cc}
                    I_{2N} & \\ & O_2 \end{array} \right),
$$
where $\mathcal{R}_1$ is obtained by adding to $\mathcal{R}$ two rows whose entries are all zeros except for \begin{eqnarray*}
&& (\mathcal{R}_1)_{2N+1,2N+1} = - (\mathcal{R}_1)_{2N+1,2N+2} =1,    \\
&& (\mathcal{R}_1)_{2N+2,2(N-s)+1} = (\mathcal{R}_1)_{2N+2,2(N-s)+2} =
\beta_{2k-s}^{(2k)}, \quad s=0,1,\dots,2k.
\end{eqnarray*}

Concerning TCII, the treatment of the boundary condition in (\ref{bcnoi}) is simpler:
it is sufficient to remove the last entry of the vector $\hat{\mathbf{v}}$ thus obtaining the
vector $\mathbf{v}$ (see (\ref{vetv}))
and the last column of the matrices $\hat{\Gamma}, \, \hat{B}, \hat{D}_0.$ More precisely, by setting
$$
\hat{\Gamma} =[ \Gamma \, | \, {\boldmath{\gamma}}_{N+1}], \quad \hat{B}= [ B
\, | \, {\boldmath{\beta}}_{N+1}], \quad \mathcal{R}_2 =  -D_2 \Gamma +D_1 B +D_0,
$$
the algebraic EVP reads:
$$\mathcal{R}_2 \,  \mathbf{v} = \lambda  \mathbf{v}.$$

\section{Numerical experiments} \label{sec4}
For the numerical simulations we considered the following potentials:
$$
V_1(r) = - \frac{2}{r}, \qquad
V_2(r) =  - \frac{2\, \alpha  \, e^{- \alpha r}}{1- e^{- \alpha r}}, \, \alpha >0, \qquad
V_3(r) = - \frac{2\,  e^{- \alpha r}}{r},\, \alpha >0,
$$
hydrogen atom, Hulth\'{e}n potential, and Yukawa potential, respectively. When $r$ is close to zero,
these three potentials  behave similarly, i.e. $V_j(r)\sim V_1(r)$ for $j=2,3.$
On the other hand,  $|V_2(r)|$ and $|V_3(r)|$ decrease faster than $|V_1(r)|$
when $r \rightarrow \infty.$ \\
For the hydrogen atom problem, the exact eigenvalues are known to be $\lambda_n = -{n^{-2}}, \, n\ge \ell+1$,
where $n$ and $\ell$ represent the radial and the angular momentum quantum
numbers, respectively,  and the corresponding eigenfunction $u_n(r)$
has exactly $\nu = n-\ell-1$ zeros in $(0,\infty).$ In the terminology of
Sturm-Liouville problems, $\lambda_n$ has therefore index $\nu.$
We solved this problem with various values of $\ell$ by applying the $(2k)$-step
scheme of order $p=2k$
described in the previous section with different values of $k$ and different
numbers of interior meshpoints $N.$ The resulting generalized eigenvalue
problems have been solved by using the {\tt eig} routine of {\sc Matlab}.
When dealing with TCII, numerical experiments indicate that a good heuristic law
for the choice of the parameter $\xi$ is given by
\begin{equation}\label{sceltaxi}
\xi = (1.35)^p \left(\ell+1\right).
\end{equation}
There are various alternative possibilities to compress the semi-infinite
interval to a finite domain.
Any transformation of the type (ATCII),
$$
t(r) = 1 - (1+r)^{-\beta}, \quad \beta > 0, \quad r \in [0,\infty),
$$
reduces $[0,\infty)$ to $[0,1)$. To see how this transformation performs in the
context of EVPs, we used ATCII with $\beta = \frac{1}{2}$. For the respective
analysis, we refer the reader to \cite{acetoetal2013}.\\

\begin{figure}
\begin{center}
      \includegraphics[width=12.5cm,height=9.5cm]{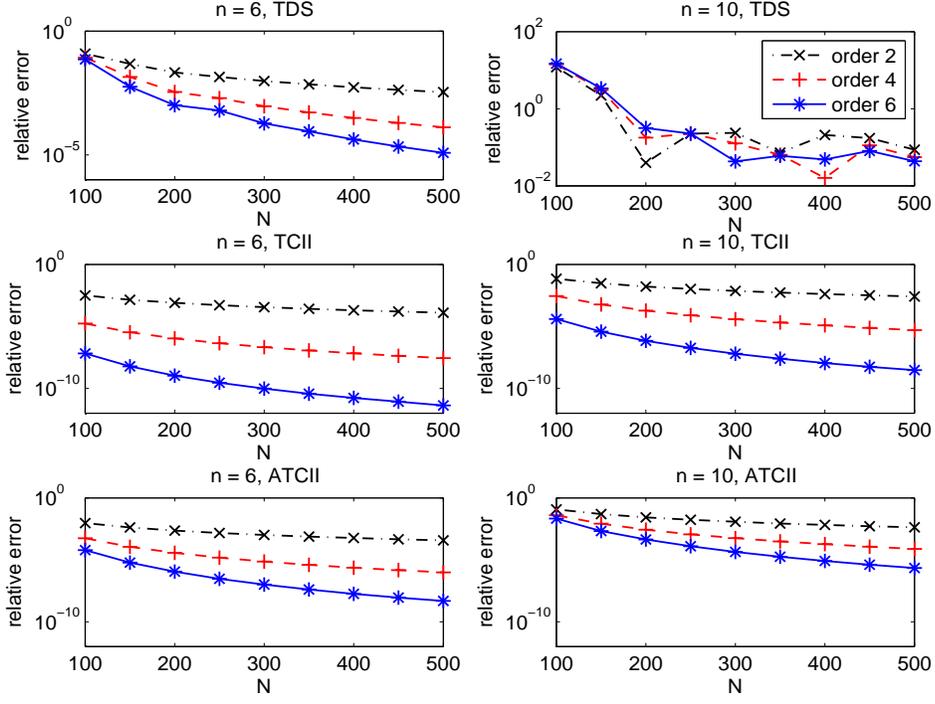}
\end{center}
\caption{\label{fig1} Hydrogen atom equation, $\ell=3$: relative errors in
the eigenvalues calculated using TDS, TCII with (\ref{sceltaxi}), and ATCII.}
\end{figure}

In Figure~\ref{fig1}, we plotted the relative errors in the eigenvalues
$\lambda_6$ and $\lambda_{10}$ of the hydrogen atom problem with $\ell=3$ versus
$N.$ In particular, the plots at the top of the picture refer to TDS, those
in the center to TCII and (\ref{sceltaxi}), and the bottom ones to ATCII. We can
see that when the radial quantum number $n$ increases, the results obtained with
TDS are not satisfactory, even for higher order methods. For TCII, we obtain
good results using already a second order method. They can be further improved
when we increase the order of the scheme. By virtue of these results and taking
into account that for a fixed $N$ the size of the generalized eigenvalue problem
corresponding to TDS is approximately twice as large as the one corresponding to
TCII, we do not include TDS in the sequel. Also, the accuracy obtained using
TCII is considerably better than the accuracy of ATCII. \\

\begin{figure}
\begin{center}
      \includegraphics[width=12.5cm,height=7cm]{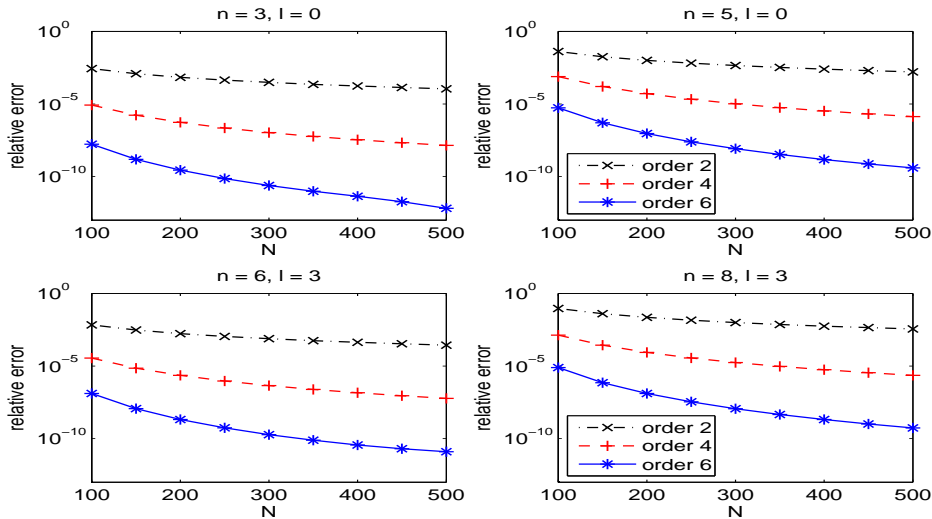}
\end{center}
\caption{\label{fig2} Hulth\'{e}n potential, $\ell=0,3$, $\alpha=0.02$:
relative errors in the eigenvalues calculated using TCII with (\ref{sceltaxi}).
}
\end{figure}

\begin{figure}
\begin{center}
      \includegraphics[width=12.5cm,height=7cm]{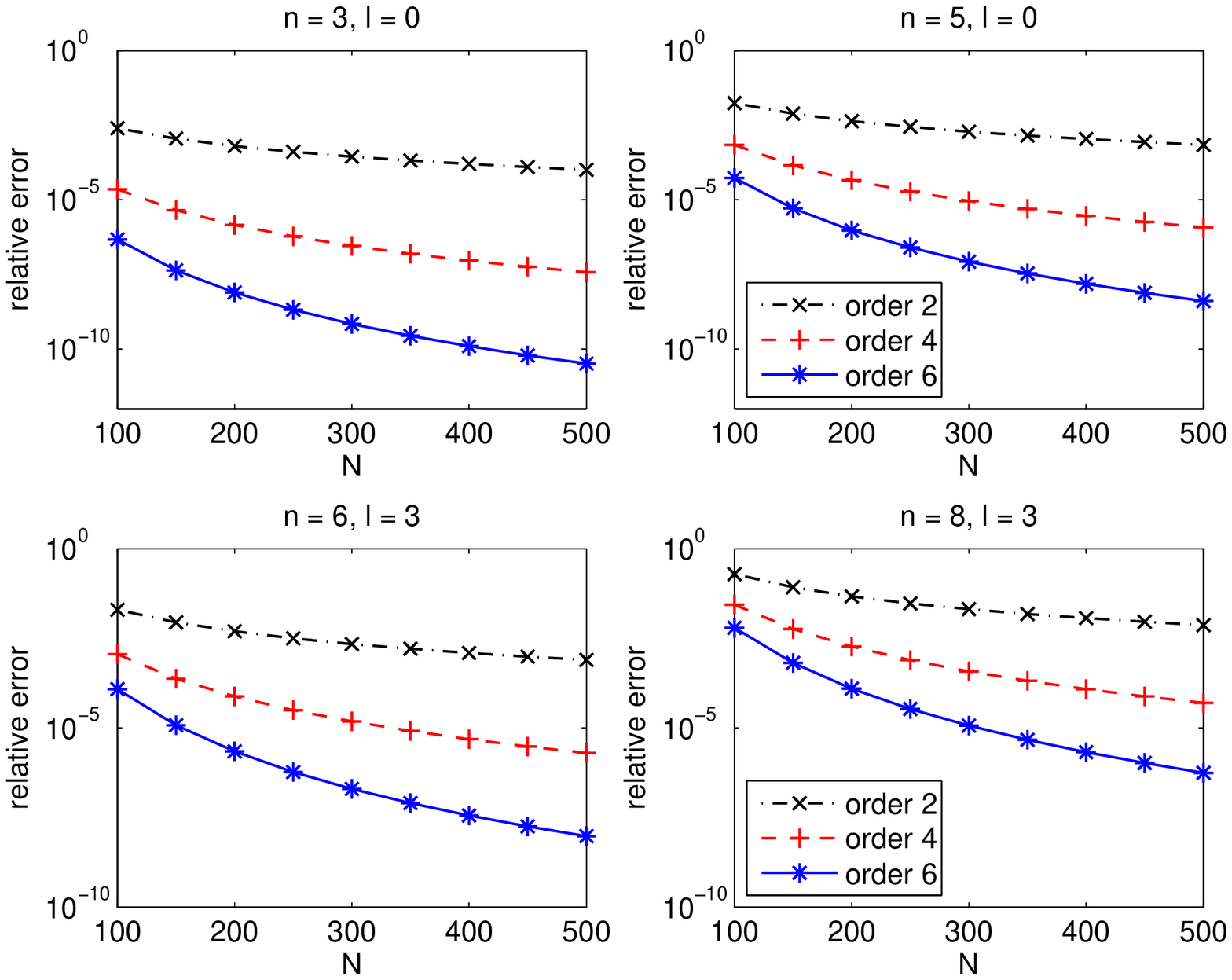}
\end{center}
\caption{Hulth\'{e}n potential, $\ell=0,3$, $\alpha=0.02$:
relative errors in the eigenvalues calculated using ATCII.} \label{ATCIIHulthen}
\end{figure}

Let us now consider the Hulth\'{e}n and the Yukawa potentials. The parameter $\alpha$ occurring in
their definition is called {\em screening parameter} and it is known that the number of
eigenvalues in the point spectrum of the corresponding problems varies with $\alpha$ \cite{Roy}.
Concerning the exact eigenvalues, these are known in closed form only for the Hulth\'en problem with
$\ell=0.$ In all other cases, in order to evaluate the performance of our schemes, we calculated the
reference eigenvalues using the method of order $8$ with $N=1500.$ As an example, in Figure~\ref{fig2},
the relative errors in the Hulth\'{e}n eigenvalue approximations for $\ell=0,3$ and $\alpha=0.02$ are shown.
Observe that both plots on the left refer to the eigenvalues of index $\nu =
n-\ell-1=2$ while the plots on the right to those of index $\nu=4.$\\ The
related data for ATCII can be found in Figure~\ref{ATCIIHulthen}.\\

In Table~\ref{tab1}, the eigenvalue approximations computed with
TCII and (\ref{sceltaxi}) using the method of order $p=8$ for the Yukawa potential have been listed and
compared to those provided by \cite{Roy}.\\

\begin{table}
\caption{Yukawa potential: eigenvalues calculated using the method of order
$p=8$ as compared the those listed in Table~5 in \cite{Roy}.}\label{tab1}
\begin{center}
\begin{small}
\begin{tabular}{|c|c|c|c|c|}
\hline
 \multicolumn{5}{|c|}{$n=9$}  \\
\hline
$\ell$ & $\alpha$ &  $\lambda/2, \quad N=200$ & $\lambda/2, \quad N=1500$  & \cite{Roy}  \\
\hline
0 & 0.010 &   -0.0005858266584  &  -0.0005858247613  &-0.0005858247612\\
1 & 0.010 &   -0.0005665076452  &  -0.0005665076262  &-0.0005665076261 \\
2 & 0.010 &   -0.0005276644219  &  -0.0005276644203  &-0.0005276644203 \\
3 & 0.010 &   -0.0004688490639  &  -0.0004688490636  &-0.0004688490636 \\
4 & 0.010 &   -0.0003893108560  &  -0.0003893108559  &-0.0003893108558 \\
5 & 0.010 &   -0.0002878564558 &   -0.0002878564558  &-0.0002878564558 \\
6 & 0.005 &   -0.0022606077423 &   -0.0022606077423  &-0.0022606077422 \\
7 & 0.005 &   -0.0021997976659 &   -0.0021997976659  &-0.0021997976659 \\
8 & 0.005 &   -0.0021291265596 &   -0.0021291265596  & -0.0021291265596\\
\hline
 \multicolumn{5}{|c|}{$n=10$}  \\
\hline
$\ell$ & $\alpha$ &  $\lambda/2, \quad N=200$ & $\lambda/2, \quad N=1500$  & \cite{Roy}  \\
\hline
0 & 0.005 &  -0.0015083751962 & -0.0015083559308 &-0.0015083559307\\
1 & 0.005 &  -0.0015009237055 & -0.0015009235029 &-0.0015009235029\\
2 & 0.005 &  -0.0014860116411 & -0.0014860116241 & -0.0014860116240\\
3 & 0.005 &  -0.0014635239308 & -0.0014635239276 & -0.0014635239275\\
4 & 0.005 &  -0.0014333097815 & -0.0014333097805 &-0.0014333097805\\
5 & 0.005 &  -0.0013951561297 & -0.0013951561294 & -0.0013951561294\\
6 & 0.005 &  -0.0013487749861 & -0.0013487749860 & -0.0013487719860\\
7 & 0.005 &  -0.0012937846260 & -0.0012937846260 & -0.0012937846259\\
8 & 0.005 &  -0.0012296811836 & -0.0012296811836 & -0.0012296811835 \\
\hline
\end{tabular}
\end{small}
\end{center}
\end{table}

\section{Conclusions}
In this paper we studied the numerical solution of the eigenvalue problems for singular Schr\"odinger equation
posed on a semi-infinite interval
\[ - u''(r)+\left(  \frac{\ell (\ell+1)}{r^2} + V(r) \right) u(r) = \lambda
u(r), \quad u(0) = u(\infty)=0. \]
Our aim was to propose a transformation reducing the infinite domain to the finite interval $(0,1]$ and
then discretize the resulting ODE using finite difference schemes. Finally,
the generalized algebraic eigenvalue problem was solved using the eigenvalue {\sc Matlab} routine.
Three transformations have been used:
\begin{enumerate}
\item[TDS:] Here, the interval $(0, \infty)$ is split into two parts, $(0,
\infty) = (0,1] \cup [1,\infty),$
and the second interval is transformed to $(0,1]$ using $t(r) := 1/r.$
This transformation has two disadvantages: the number of equations is doubled
which is not so critical since the original problem is scalar, but also a
singularity of the fist kind in the original problem changes to an essential
singularity in the transformed equations.
The latter singularity is considerably more difficult to handle numerically.
    \item[TCII:] With the transformation $t(r):=r/(r + \xi)$ and a suitably
chosen $\xi,$ the semi-infinite interval is compressed to $(0,1)$; the
dimension of the problem and the type of the singularity do not change.
\item[ATCII:] Analogous compression is also done using $t(r):=1 - 1/\sqrt{1+r}$.
\end{enumerate}
We could show that the transformed problems are well-posed and discussed the smoothness of their solutions.
Moreover, it turns out that the approach based on TCII outperforms the other two, and therefore,
it could be recommended to be used in similar situations.

\section*{Acknowledgements}
The authors wish to thank Pierluigi Amodio and Giuseppina Settanni for providing
the software for setting up the finite difference schemes.

\end{document}